\documentclass[11pt,a4paper]{amsart}

\usepackage[margin=27mm]{geometry}
\usepackage[T1]{fontenc}
\usepackage{lmodern,amssymb,microtype}
\usepackage[hidelinks]{hyperref}
\hypersetup{pdftitle={Equality in the direction bound over prime fields},
  pdfauthor={P\'eter Sziklai}}

\newtheorem{theorem}{Theorem}
\newtheorem{lemma}[theorem]{Lemma}
\newtheorem{proposition}[theorem]{Proposition}
\theoremstyle{remark}

\newcommand{\F}{\mathbb F_p}
\newcommand{\AG}{\operatorname{AG}}
\newcommand{\LL}{\mathcal L}
\newcommand{\mult}{\operatorname{mult}}
\newcommand{\tot}{\mathrm{tot}}

\title{Equality in the direction bound over prime fields}
\author{Peter Sziklai}
\address{Institute of Mathematics, ELTE E\"otv\"os Lor\'and University,
Budapest, Hungary}
\email{peter.sziklai@ttk.elte.hu}
\date{}

\begin{document}
\begin{abstract}
Let $p$ be prime. We prove that a noncollinear set of $2m+1\le p$
points in $\AG(2,p)$ determines exactly $m+2$ directions if and only if,
up to an affine transformation, it consists of the origin and a coset
of a multiplicative subgroup of order $m$ on each coordinate axis.
Thus equality is possible precisely when $m\mid p-1$. This extends
the Lov\'asz--Schrijver classification for $2m+1=p$.
\end{abstract}
\maketitle

\section{Introduction}

For $U\subseteq\AG(2,p)=\F^2$, let $D(U)\subseteq\F\cup\{\infty\}$
be the set of slopes of lines joining two distinct points of $U$.
A theorem of Sz\H{o}nyi \cite{Szonyi1996} (see also
\cite[Theorem~5.2]{Szonyi}), obtained independently by Blokhuis (see \cite[Theorem~4]{FST}), states that every noncollinear
set with $|U|\le p$ satisfies
\[
 |D(U)|\ge \frac{|U|+3}{2}.
\]
For $|U|=p$, the bound is due to R\'edei and Megyesi \cite{Redei},
and Lov\'asz and Schrijver \cite{LS} classified the equality case. We determine all equality cases of (obviously) odd cardinality.

\begin{theorem}\label{thm:main}
Let $p$ be prime, let $m\ge1$, and let $U\subseteq\AG(2,p)$ be
noncollinear with $|U|=2m+1\le p$. Then $|D(U)|=m+2$ if and only if
$m\mid p-1$ and an affine transformation takes $U$ to
\[
 \{(0,0)\}\ \cup\ (H\times\{0\})\ \cup\ (\{0\}\times H),
\]
where $H$ is the subgroup of order $m$ in $\F^\times$.
\end{theorem}

The configurations in Theorem~\ref{thm:main} are the Megyesi-type examples
of Blokhuis; see \cite[following Theorem~2.13]{DeBeule}.
Before independent scaling of the two axes, the two nonzero coordinate
sets are arbitrary cosets of $H$. In particular, for fixed $p$ and $m$
there is exactly one affine equivalence class of equality examples
whenever $m\mid p-1$. At $m=(p-1)/2$, the subgroup $H$ is the set of
nonzero squares, giving the Lov\'asz--Schrijver configuration.
Deleting any point from such a configuration with $m\ge2$ gives a
$2m$-point set attaining the rounded bound $m+2$; we do not treat the
even case.

The proof starts with equality in the R\'edei polynomial argument.
A second derivative supplies nonnegative integer weights whose sum on
each secant is its number of points minus two. We then prove concurrence
of the lines containing the centroids of parallel sections, adapting
the recent nice centroid-line method in Csajb\'ok \cite[Lemma~3.2]{Csajbok}.
The weighted centre identifies the concurrence point even when $|U|=p$.
A counting argument forces the positive-weight points onto a line,
and the subgroup description follows. All ingredients used in the
proof are established below.

The case $m=1$ is the affine normalization of a noncollinear triangle.
Henceforth assume $m\ge2$ and $|D(U)|=m+2$.
A \emph{secant} is a line meeting $U$ in at least two points.
Write $\LL$ for the set of secants, $k_\ell=|U\cap\ell|$, and
$t(P)$ for the number of secants through $P\in U$.
Counts and multiplicities are integers; polynomial calculations and
centroids are over $\F$. In line equations, $(\alpha,\beta)$ denote affine
coordinate variables; the constant $e$ in $\beta=d\alpha+e$ is its intercept.

\section{Equality in the R\'edei argument}

Choose coordinates with $\infty\in D(U)$, and set
\[
 R(X,Y)=\prod_{(a,b)\in U}(X+aY-b),\qquad
 \Delta(Y)=\prod_{d\in D(U)\setminus\{\infty\}}(Y-d).
\]
Thus $\deg\Delta=m+1$. Regard \(R(X,Y)\) as a monic polynomial of degree \(2m+1\) in \(X\) over \(\mathbb F_p[Y]\). Division by $R$ gives
\begin{equation}\label{eq:division}
 RQ=X^p-X+g,\qquad \deg_Xg<2m+1,\qquad \deg_{\tot}g\le p.
\end{equation}
The total-degree bound follows directly during monic division,
since the coefficient of $X^{2m+1-i}$ in $R$ has degree at most $i$.
For every \emph{nondetermined finite direction} $y$, the polynomial $R(X,y)$
has distinct roots in $\F$, so it divides $X^p-X$. Thus all coefficients
of $g$ vanish at the $p-m-1$ such values of $y$, and
\begin{equation}\label{eq:S}
 g=WS,\qquad
 W(Y)=\frac{Y^p-Y}{\Delta(Y)},\qquad \deg_{\tot}S\le m+1.
\end{equation}

Fix a \emph{finite determined direction} $d$, and write $g_d(X)=g(X,d)$ and
$R_d(X)=R(X,d)$. If $e$ is a root of $R_d$ of multiplicity $k$, then
$g_d(e)=0$, while differentiating \eqref{eq:division} shows that
$(X-e)^{k-1}$ divides $g_d'-1$. Consequently
$R_d\mid g_d(g_d'-1)$. Neither factor is zero: $g_d=0$ would make
$R_d$ squarefree, whereas $g_d'=1$ would give $g_d=X+\gamma$,
since $\deg g_d<p$, and force $R_d\mid X^p+\gamma$, making $U$
collinear. As $2m+1=2(m+1)-1$ and $\deg g_d\le m+1$, degree comparison gives
\begin{equation}\label{eq:univariate}
 g_d(g_d'-1)=(m+1) a_d^2R_d,\qquad \deg g_d=m+1,
\end{equation}
where $a_d\ne0$ is the leading coefficient of $g_d$.

Every root of $g_d$ lies in $\F$ by \eqref{eq:univariate}.
Moreover, $g_d$ is squarefree. Indeed, a multiple root would also be a
multiple root of $R_d$, whereas the derivative of $X^p-X+g_d$ there
would be $-1$, contradicting \eqref{eq:division}.
Conversely, every root of $R_d$ is a root of $g_d$.
It follows that exactly $m+1$ lines in direction $d$ meet $U$.
Changing coordinates proves the same assertion in the vertical
direction.

For each determined direction, summing $k_\ell-1$ over its secants
therefore gives $(2m+1)-(m+1)=m$. Counting pairs of points as well, we obtain
\begin{equation}\label{eq:counts}
 \sum_{\ell\in\LL}(k_\ell-1)=m(m+2),\qquad
 \sum_{\ell\in\LL}\binom{k_\ell}{2}=m(2m+1).
\end{equation}
In particular,
\begin{equation}\label{eq:moment}
 \sum_{\ell\in\LL}(k_\ell-1)(k_\ell-2)=2m(m-1).
\end{equation}
Also $k_\ell\le m+1$: a point of $U$ outside $\ell$ gives $k_\ell$
distinct joining directions, in addition to the direction of $\ell$.

The coefficient of $X^{m+1}$ in $S$ is constant since
$\deg_{\tot}S\le m+1$, and is nonzero by \eqref{eq:univariate} and
$W(d)\ne0$ for any finite determined $d$; denote it by $c$.
Substituting $X=b-aY$ in \eqref{eq:division} gives
\begin{equation}\label{eq:restriction}
 S(b-aY,Y)=a\Delta(Y)\qquad ((a,b)\in U).
\end{equation}
On the line $X=b-aY$ we have $S=a\Delta$ and, by differentiating
\eqref{eq:restriction}, $S_Y=a(S_X+\Delta')$.
Hence every linear factor of $R$ divides $S(S_X+\Delta')-\Delta S_Y$. Comparing total degrees and
the coefficient of $X^{2m+1}$ yields
\begin{equation}\label{eq:identity}
 S(S_X+\Delta')-\Delta S_Y=(m+1)c^2R.
\end{equation}

\section{Integer weights}

Homogenize $S$ and $\Delta$ to degree $m+1$ in $X,Y,Z$, writing
$\overline S$ and $\overline\Delta$, and put
\[
 \Phi=\overline S_{XX}.
\]
This is a nonzero homogeneous polynomial of degree $m-1$, with
coefficient of $X^{m-1}$ being $m(m+1)c$. In the projective dual plane,
$P=(a,b)$ corresponds to the line
\[
 L_P=X+aY-bZ=0.
\]
Let $w_P$ be the exact multiplicity of the factor $L_P$ in $\Phi$,
allowing $w_P=0$.

\begin{lemma}\label{lem:local}
For every $k$-secant $\ell$, with dual point $\ell^*$,
\[
 \mult_{\ell^*}\Phi\ge k-2,\qquad
 \sum_{P\in U\cap\ell}w_P\le k-2.
\]
\end{lemma}

\begin{proof}
Suppose first that $\ell$ has equation $\beta=d\alpha+e$ in the variables $\alpha, \beta$. Write
$\overline\Delta=(Y-dZ)\delta(Y,Z)$. By \eqref{eq:restriction},
\[
 \overline S+(X-eZ)\delta(Y,Z)
\]
is divisible by $\prod_{P\in U\cap\ell}L_P$, and hence has multiplicity
at least $k$ at their common point $\ell^*$. Its second $X$ derivative
is $\Phi$, proving the first inequality. For a vertical line $\alpha=a_0$,
use $\overline S-a_0\overline\Delta$ instead.

For the second inequality, fix $P=(a_0,b_0)$ with $w_P=w>0$ and
translate $P$ to the origin. The corresponding polynomial is
$S(X-a_0Y+b_0,Y)-a_0\Delta(Y)$; it satisfies the same identities and
preserves the multiplicity in the second derivative.
Thus we may assume $P=(0,0)$.
By \eqref{eq:restriction}, $X\mid S$. Since $X^w\mid S_{XX}$ and
$m+1<p$, the coefficients of $X^2,\ldots,X^{w+1}$ in $S$ vanish.
Equation \eqref{eq:identity} shows that the same coefficients in $R$
vanish.

Apply Newton's identities over $\F(Y)$ to the reciprocal roots of
\[
 R/X=\prod_{(a,b)\in U\setminus\{P\}}(X+aY-b).
\]
The vanishing coefficients imply
\begin{equation}\label{eq:reciprocal}
 \sum_{(a,b)\ne P}(aY-b)^{-j}=0\qquad(1\le j\le w).
\end{equation}
Group the terms by lines through $P$. For $a\ne0$, the term
$(aY-b)^{-j}$ has its pole at $Y=b/a$, the slope of the line through
$P$ and $(a,b)$. Distinct nonvertical lines give distinct poles in $Y$,
so their coefficients vanish separately;
the vertical group is the remaining constant term.
On any $k$-secant $\ell$ through $P$, choose an affine coordinate with $P=0$
and other point-coordinates $u_1,\ldots,u_{k-1}$. We obtain
$\sum_i u_i^{-j}=0$ for $1\le j\le w$.
Newton's identities now show that
\[
 f_\ell(T)=T\prod_{i=1}^{k-1}(T-u_i)
\]
has zero coefficients at $T^2,\ldots,T^{w+1}$.
Thus $f_\ell''$ has a root of multiplicity at least $w$ at $P$.
All the integers divided by here are less than $p$.
Translation of the affine coordinate preserves root multiplicities
of the second derivative. Thus, in one common coordinate on $\ell$,
repeating the argument gives a root of $f_\ell''$ of multiplicity at
least $w_Q$ at each positive-weight point $Q\in U\cap\ell$. Since $f_\ell''$ is nonzero of degree $k-2<p$,
the sum of these multiplicities is at most $k-2$.
\end{proof}

\begin{proposition}\label{prop:weights}
The polynomial $\Phi$ factors completely as
\[
 \Phi=m(m+1)c\prod_{P\in U}L_P^{w_P}.
\]
The nonnegative integer weights satisfy
\begin{align}
 \sum_{P\in U}w_P&=m-1,\label{eq:weightsum}\\
 \sum_{P\in U\cap\ell}w_P&=k_\ell-2
       &&(\ell\in\LL),\label{eq:lineweight}\\
 (t(P)-1)(w_P+1)&=m
       &&(P\in U).\label{eq:divisor}
\end{align}
In particular, the weights are invariant under affine transformations.
\end{proposition}

\begin{proof}
Write $\Phi=(\prod_P L_P^{w_P})T$, where no $L_P$ divides $T$, and set
\[
 e=\deg T=m-1-\sum_Pw_P,\qquad
 \delta_\ell=k_\ell-2-\sum_{P\in U\cap\ell}w_P.
\]
Lemma~\ref{lem:local} gives $\delta_\ell\ge0$ and
$\mult_{\ell^*}T\ge\delta_\ell$. The restriction of $T$ to each $L_P$
is a nonzero homogeneous polynomial of degree $e$ in two variables.
Its zeros at the dual points of secants through $P$ therefore give
$\sum_{\ell\ni P}\delta_\ell\le e$.
Summing over $P$ yields
\begin{equation}\label{eq:bezout}
 \sum_{\ell\in\LL}k_\ell\delta_\ell\le(2m+1)e.
\end{equation}
On the other hand, \eqref{eq:moment} and
$\sum_{\ell\ni P}(k_\ell-1)=2m$ imply
\begin{equation}\label{eq:deficit}
 \sum_{\ell\in\LL}(k_\ell-1)\delta_\ell
 =2m(m-1)-2m\sum_Pw_P=2me.
\end{equation}
Subtracting \eqref{eq:deficit} from \eqref{eq:bezout} gives
$\sum_\ell\delta_\ell\le e$. Since $k_\ell-1\le m$,
\[
 2me\le m\sum_\ell\delta_\ell\le me.
\]
Hence $e=0$, and \eqref{eq:deficit} forces every $\delta_\ell=0$.
This proves the factorization, \eqref{eq:weightsum}, and
\eqref{eq:lineweight}.
Finally, summing \eqref{eq:lineweight} over the secants through $P$
gives
\[
 m-1+(t(P)-1)w_P=2m-t(P),
\]
which is \eqref{eq:divisor}. This formula determines $w_P$ from the
affine-invariant number $t(P)$.
\end{proof}

For the configurations in Theorem~\ref{thm:main}, all weight is concentrated
at the origin: it lies on exactly two secants, so its weight is $m-1$
by \eqref{eq:divisor}.

\section{Centroid lines and the classification}

Every nonempty line intersection has at most $m+1<p$ points, so its
centroid is defined. Introduce the weighted centre
\begin{equation}\label{eq:centre}
 C=\frac1{m-1}\sum_{P\in U}w_PP.
\end{equation}
By Proposition~\ref{prop:weights}, this point transforms covariantly
under affine maps. The denominator $m-1$ is nonzero even when $|U|=p$.

\begin{lemma}\label{lem:centroid}
For each $d\in D(U)$, the centroids of the nonempty intersections of
$U$ with lines of direction $d$ lie on a line $M_d$ through $C$.
The direction of $M_d$ differs from $d$.
\end{lemma}

\begin{proof}
First let $d$ be finite, in the coordinates of Section~2. Write
\[
 S(X,Y)=cX^{m+1}+s(Y)X^{m}+\cdots,\qquad \deg s\le1.
\]
Comparing the first two $X$ coefficients in the factorization of
$\overline S_{XX}$ gives
\[
 \sum_Pw_P(a_PY-b_P)=\frac{m-1}{(m+1)c}s(Y).
\]
Thus, if $C=(u,v)$, then
\begin{equation}\label{eq:centrecoeff}
 \frac{s(Y)}{(m+1)c}=uY-v.
\end{equation}

Set
\[
 \delta(Y)=\frac{\Delta(Y)}{Y-d},\qquad
 \eta=\frac{\delta'(d)}{\delta(d)}
      =\frac{\Delta''(d)}{2\Delta'(d)},\qquad
 B_d(X)=S_X(X,d)+\Delta'(d).
\]
Here $2$ is invertible since $p\ge2m+1\ge5$.
Differentiating $Y^p-Y=W\Delta$ at $Y=d$ gives $W(d)\Delta'(d)=-1$.
Thus \eqref{eq:univariate} and the squarefreeness
of $S(X,d)$ show that the roots of $B_d$ are precisely the secant
intercepts in direction $d$, with multiplicity $k_\ell-1$.
Also $\deg B_d=m$.

Consider an occupied line $\beta=d\alpha+e$, whose $k$ points have abscissas
$a_1,\ldots,a_k$. As in Lemma~\ref{lem:local},
\begin{equation}\label{eq:localfactor}
 S(X,Y)+(X-e)\delta(Y)
 =\prod_{i=1}^k\bigl(X+a_i(Y-d)-e\bigr)V(X,Y).
\end{equation}
Put $N_d(X)=S_Y(X,d)-\eta S(X,d)$.
Differentiating \eqref{eq:localfactor} with respect to $X$, and applying
$\partial_Y-\eta$ to it, gives respectively
\begin{align*}
 B_d(X)&\equiv kV(e,d)(X-e)^{k-1}\pmod{(X-e)^k},\\
 N_d(X)&\equiv\left(\sum_i a_i\right)V(e,d)(X-e)^{k-1}
                         \pmod{(X-e)^k}.
\end{align*}
Here $V(e,d)\ne0$, since the multiplicity in $B_d$ is exactly $k-1$;
this includes $k=1$, when $B_d(e)\ne0$.
It follows that $B_d\mid N_d$ and
\[
 A_d(X):=\frac{N_d(X)}{B_d(X)}
 \quad\hbox{satisfies}\quad
 \deg A_d\le1,\qquad A_d(e)=\frac1k\sum_i a_i.
\]
Comparison of the two highest coefficients, using
\eqref{eq:centrecoeff}, gives
\[
 A_d(X)=u-\frac{\eta}{m+1}(X+ud-v).
\]
If $(\bar a,\bar b)$ is the centroid on the line with intercept $e$,
then $e=\bar b-d\bar a$. Hence
\[
 (m+1)(\bar a-u)+\eta\bigl((\bar b-v)-d(\bar a-u)\bigr)=0.
\]
This is a line through $C$ whose direction differs from $d$, since
substitution of the direction vector $(1,d)$ gives $m+1\ne0$.
An affine change of coordinates treats the vertical direction as well;
the covariance of \eqref{eq:centre} ensures the same centre $C$.
\end{proof}

\begin{lemma}\label{lem:collinear}
The set $U^+=\{P\in U:w_P>0\}$ is collinear.
\end{lemma}

\begin{proof}
The set $U^+$ is nonempty by \eqref{eq:weightsum}.
For $P\in U$, let $E(P)$ be the set of determined directions in which
the line through $P$ meets $U$ only at $P$.
Equation \eqref{eq:divisor} gives
\begin{equation}\label{eq:singleton}
 |E(P)|=m+2-t(P)=m+1-\frac{m}{w_P+1}.
\end{equation}
For $d\in E(P)$, the point $P$ is itself a centroid, so
$P\in M_d$ by Lemma~\ref{lem:centroid}.
If $P,Q\in U^+$ and $C,P,Q$ are noncollinear, then
$E(P)\cap E(Q)=\varnothing$. Thus \eqref{eq:singleton} implies
\[
 2m+2-\frac{m}{w_P+1}-\frac{m}{w_Q+1}\le m+2,
 \qquad\text{or}\qquad
 \frac1{w_P+1}+\frac1{w_Q+1}\ge1.
\]
Since $w_P+1,w_Q+1\ge2$, we obtain $w_P=w_Q=1$.

Suppose, for a contradiction, that $U^+$ is not collinear.
For every $P\in U^+\setminus\{C\}$ there is a $Q\in U^+$ outside the
line $CP$. The preceding argument shows that $w_P=1$.
Consequently $m$ is even and
$|E(P)|=(m+2)/2$ for every such $P$.
Points on different lines through $C$ have disjoint sets $E(P)$.
There can therefore be at most two such lines containing points of
$U^+\setminus\{C\}$; noncollinearity forces exactly two, say $L_1,L_2$.

Put $B_i=(U^+\setminus\{C\})\cap L_i$ and $b_i=|B_i|>0$.
Let $\omega=w_C$ if $C\in U$, and $\omega=0$ otherwise.
Then
\begin{equation}\label{eq:Bmass}
 \omega+b_1+b_2=m-1.
\end{equation}
The centre identity gives
\[
 \sum_{P\in B_1}(P-C)+\sum_{P\in B_2}(P-C)=0.
\]
The two sums belong to distinct one-dimensional subspaces, so both
vanish. In particular, $b_i\ge2$. Thus $L_i$ is a secant, and
\eqref{eq:lineweight} gives
\[
 |U\cap L_i|=\omega+b_i+2.
\]
A point $P\in B_1$ lies outside $L_2$, and its joins to the points of
$U\cap L_2$ are distinct secants. Since $w_P=1$,
\eqref{eq:divisor} yields
$\omega+b_2+2\le t(P)=m/2+1$.
Interchanging the indices and adding gives
$2\omega+b_1+b_2\le m-2$, contradicting \eqref{eq:Bmass}.
\end{proof}

\begin{proof}[Proof of Theorem~\ref{thm:main}]
For $m\ge2$, choose a secant containing all of $U^+$; if $U^+$ is a
singleton, any secant through it will do.
Equations \eqref{eq:weightsum} and \eqref{eq:lineweight} show that
this line contains $m+1$ points of $U$.

Put this line at $y=0$, with abscissa set $A$, $|A|=m+1$.
Parametrize the nonhorizontal determined directions by horizontal displacement
divided by vertical displacement, and denote their set by $T$.
Then $|T|=m+1$, and each point $(u,v)\in U$ outside $y=0$ gives
\[
 (u-A)/v=T.
\]
Taking arithmetic means in $\F$ gives
$u-v\,\operatorname{mean}(T)=\operatorname{mean}(A)$. This is an equation for $(u,v)$;
hence all $m$ outside points lie on one nonhorizontal line $M$.
Since $m\ge2$, its direction is determined and belongs to $T$.
The displayed equality then implies that the intersection point of
$M$ with $y=0$ belongs to $U$.

Send these two lines to the coordinate axes, and their intersection
to the origin. Their nonzero coordinate sets $A_0,B_0$ each have
size $m$. The cross-axis directions give
\[
 |B_0/A_0|=m.
\]
For fixed $a_0\in A_0$, all sets $B_0/a$ with $a\in A_0$ are equal,
since each has size $m$ and is contained in $B_0/A_0$.
Hence every $a/a_0$ stabilizes $B_0$ multiplicatively.
Its stabilizer $H\le\F^\times$ has at least $m$ elements, and at most
$m$ because its orbits on $B_0$ have size $|H|$.
Therefore $|H|=m$, both $A_0$ and $B_0$ are cosets of $H$, and
$m\mid p-1$. Independent scaling of the axes gives the stated form.

Conversely, for a subgroup $H$ of order $m$, the stated set has
$2m+1$ points and direction set $\{0,\infty\}\cup(-H)$.
It therefore determines exactly $m+2$ directions.
The case $m=1$, already noted in the introduction, completes the proof.
\end{proof}

\bigskip
\noindent{\bf Use of AI.}
This research was assisted by AI tools, in parts of mathematical exploration, computations
and drafting. The author takes full responsibility for
the manuscript.

\bigskip
\noindent{\bf Acknowledgements.} Project no. 153080 has been implemented with the support provided
by the Ministry of Culture and Innovation of Hungary from the National Research, Development and Innovation Fund, financed under the ADVANCED 25 funding schemes.

\end{document}